\documentclass[12pt]{article}
\usepackage{latexsym,amssymb,amsmath}
\usepackage{color}
\usepackage{cite}
\usepackage{setspace}
\usepackage{enumerate}

\newtheorem{Theorem}{\bf Theorem}[section]
\newtheorem{Lemma}{\bf Lemma}[section]
\newtheorem{Proposition}{\bf Proposition}[section]
\newtheorem{Corollary}{\bf Corollary}[section]
\newtheorem{Remark}{\bf Remark}[section]
\newtheorem{Example}{\bf Example}[section]
\newtheorem{Definition}{\bf Definition}[section]

\numberwithin{equation}{section}

\def\XXint#1#2#3{{\setbox0=\hbox{$#1{#2#3}{\int}$}
\vcenter{\hbox{$#2#3$}}\kern-.5\wd0}}

\makeatletter
\newcommand{\subjclass}[2][2010]{%
  \let\@oldtitle\@title%
  \gdef\@title{\@oldtitle\footnotetext{#1 \emph{AMS subject classification.} #2}}%
}
\newcommand{\keywords}[1]{%
  \let\@@oldtitle\@title%
  \gdef\@title{\@@oldtitle\footnotetext{\emph{Key words and phrases.} #1.}}%
}
\makeatother

\begin{document}
  
  \title{Existence of solutions to fractional semilinear parabolic equations in Besov-Morrey spaces}
\author{Erbol Zhanpeisov\thanks{{E-mail: erbol.zhanpeisov@oist.jp}}\\
Okinawa Institute of Science and Technology\\
1919-1 Tancha, Onna-son, Kunigami-gun\\
 Okinawa, Japan 904-0495}
\date{}
\maketitle
\vspace{-10mm}
\begin{abstract}
In this paper, we establish the existence of solutions to fractional semilinear parabolic equations in Besov-Morrey spaces for a large class of initial data including distributions other than Radon measures. We also obtain  sufficient conditions for the existence of solutions to viscous Hamilton-Jacobi equations.
\end{abstract}


\section{Introduction and main results}
Consider a semilinear parabolic equation 
\begin{equation}\label{eq:1.1}
\left\{
\begin{aligned}
\partial_t u  + (-\Delta)^{\frac{\theta}{2}} u = |u|^{\gamma-1}u, \quad & x\in {\bf R}^N, \ t\in (0,T), \\
u(x,0) = \varphi(x), \quad &x\in {\bf R}^N
\end{aligned}
\right.
\end{equation}
and a viscous Hamilton-Jacobi equation
\begin{equation}
\label{eq:1.2}
\left\{
\begin{aligned}
\partial_t u  + (-\Delta)^{\frac{\theta}{2}} u = |\nabla u|^{\gamma}, \quad & x\in {\bf R}^N, \ t\in (0,T), \\
u(x,0) = \varphi(x), \quad &x\in {\bf R}^N, 
\end{aligned}
\right.
\end{equation}
where $\gamma>1$, $N\geq 1$, $T>0$ and $\theta >0$ (resp. $\theta>1$)
for problem~\eqref{eq:1.1} (resp. problem~\eqref{eq:1.2}). The purpose of this paper is to obtain sufficient conditions for the existence of solutions to the Cauchy problem~\eqref{eq:1.1} and \eqref{eq:1.2} for a large class of  initial data by introducing inhomogeneous Besov-Morrey spaces. This enables us to take distributions other than Radon measures as initial data. 

Let us consider the Cauchy problem for the semilinear parabolic equation~\eqref{eq:1.1} with $\theta>0$ and $\gamma>1$. The solvability of problem~\eqref{eq:1.1} has been studied in many papers, see e.g., \cite{BP, C2, FL, G, R, Wu, GP, GG, HI1, HI2, HIT, IKO1, IKO2, KY, LN, P, Su, T, W}. 
(See also the monograph \cite{QS}.) Among others, Ishige, Kawakami, and Okabe \cite{IKO2} developed the arguments in \cite{IKO1} and  obtained sufficient conditions for the existence of solutions to problem~\eqref{eq:1.1} for general $\theta>0$. As corollaries of their main results, they proved the following properties: 
\begin{enumerate}[(a)]
\item Let $1<\gamma<1+\theta/N$. Then problem~\eqref{eq:1.1} possesses a local-in-time solution if $$\sup_{x\in {\bf R}^N}\|\varphi \|_{L^1(B(x,1))}<\infty;$$ 
\item Let $\gamma=1+\theta/N$. Then there exists $c>0$ such that, if 
\begin{equation*}
|\varphi (x)|\le c |x|^{-N}\left|\log \left(e+\frac{1}{|x|}\right)\right|^{-\frac{N}{\theta}-1}, \quad x\in{\bf R}^N, \\
\end{equation*}
then probolem~\eqref{eq:1.1} possesses a local-in-time solution; 
\item Let $\gamma>1+\theta/N$. Then there exists $c>0$ such that, if 
\begin{equation*}
|\varphi (x)|\le c |x|^{-\frac{\theta}{\gamma-1}}, \quad x\in{\bf R}^N, \\
\end{equation*}
then probolem~\eqref{eq:1.1} possesses a local-in-time solution. 
\end{enumerate}
In the case of either $0<\theta \le 2$ or $\theta\in \{4,6,\dots\}$, it is shown in \cite{HI1} and \cite{IKO1} that sufficient conditions in (b) and (c) are sharp. More precisely, there exists $c'>0$ such that, if 
\begin{equation*}
\varphi(x)\ge \left\{
\begin{aligned}
& c' |x|^{-N}\left|\log \left(e+\frac{1}{|x|}\right)\right|^{-\frac{N}{\theta}-1} \quad & \text{if} \quad \gamma = 1+\frac{\theta}{N}, \\
& c' |x|^{-\frac{\theta}{\gamma-1}} \quad & \text{if} \quad \gamma > 1+\frac{\theta}{N},
\end{aligned}
\quad x \in B(0,1),
\right.
\end{equation*}
then problem~\eqref{eq:1.1} possesses no local-in-time nonnegative solutions.

On the other hand, in the case of (a), distributions other than Radon measures such as the derivative of the Dirac distribution can be considered as the initial data to problem~\eqref{eq:1.1} with $\theta=2$. For instance, problem~\eqref{eq:1.1} with $\theta=2$ is well-posed in certain 
negative order inhomogeneous Besov-Morrey spaces $N^s_{p,q,r}({\bf R}^N)$, see \cite{KY} and Remark~\ref{Remark:1.1}. 
The arguments in \cite{KY} are based on delicate decay estimates of the heat kernel in inhomogeneous Besov-Morrey spaces and the power nonlinearity of 
the semilinear parabolic equation. It seems difficult to apply their arguments directly to the Cauchy problem~\eqref{eq:1.1} and problem~\eqref{eq:1.2}, in  particular, the case of fractional diffusion $\theta\not=2$ and the case of the nonlinearity  depending on $\nabla u$. 
In this paper, we develop the arguments in \cite{KY} and prove the unique existence of the solution to problem~\eqref{eq:1.1} (resp. problem~\eqref{eq:1.2}) in  inhomogeneous Besov-Morrey spaces $N^s_{p,q,r}({\bf R}^N)$ for general $\theta>0$ (resp. $\theta>1$).  This enables us to take distributions other than Radon measures as initial data and  the results in the case (a) is extended for more general initial data.

For viscous Hamilton-Jacobi equations~\eqref{eq:1.2}, the solvability  has been studied in \cite{AB, DI, KW, BSW, A}. 
Using the majorant kernel, Ishige, Kawakami, and Okabe \cite{IKO2} obtained the same results for problem~\eqref{eq:1.2} as for problem~\eqref{eq:1.1}.
That is, when $1<\gamma<1+(\theta+1)/(N+1)$, there exists a solution to problem~\eqref{eq:1.2} if the initial measure satisfies
$$\sup_{x\in {\bf R}^N}\|\varphi \|_{L^1(B(x,1))}<\infty.$$
We extend these results to more general initial data. See Remark~\ref{Remark:1.2} for more details on the relation to previous studies.

We recall the definition of local Morrey spaces and introduce inhomogeneous Besov-Morrey spaces.

\begin{Definition}[local Morrey spaces]
Let $1\leq q \leq p<\infty$. The local Morrey space $M^p_q({\bf R}^N)$ is defined to be the set of measurable functions $u$ in ${\bf R}^N$ such that 
\begin{equation*}
\|u~|M^p_q\|:= \sup_{x\in {\bf R}^N,~ 0<\rho \leq1}\rho^{\frac{N}{p}-\frac{N}{q}}\|u~|L^q(B(x,\rho))\|<\infty. 
\end{equation*}
The local measure space of the Morrey type $M^p({\bf R}^N)$ is defined as the sets of the Radon measures $\mu$ on ${\bf R}^N$ such that 

\begin{equation*}
\|\mu |M^p\|:=\sup_{x\in {\bf R}^N, 0<\rho \leq 1}\rho^{\frac{N}{p}-N}|\mu|(B(x,\rho))<\infty,
\end{equation*}
 where $|\mu|$ denotes the total variation of the measure $\mu$.
\end{Definition}%

Let ${\zeta}(t)$ be a smooth function on $[0,\infty)$ such that $0\leq \zeta(t)\leq1$, $\zeta(t)\equiv 1$ for $t\leq \frac{3}{2}$ and $\text{supp}~ \zeta \subset [0,\frac{5}{3})$. For $j\in {\bf Z}$, put $\varphi_j(\xi):=\zeta(2^{-j}|\xi|)-\zeta(2^{1-j}|\xi|)$ and $\varphi_{(0)}(\xi) :=\zeta(|\xi|)$.  Then we have $\varphi_j(\xi), ~\varphi_{(0)}(\xi)\in C_0^\infty ({\bf R}^N)$ and 
\begin{equation*}
\varphi_{(0)}(\xi )+ \sum_{j=1}^\infty \varphi_j(\xi)=1 \quad \text{for any}\quad \xi\in {\bf R}^N.
\end{equation*}

\begin{Definition}[inhomogeneous Besov-Morrey space]
Let $1\leq q\leq p<\infty$, $1\leq r\leq \infty$ and $s\in {\bf R}$. The local Besov-Morrey space is defined as the sets of distributions $u\in \mathcal{S}'({\bf R}^N)$ such that $\mathcal{F}^{-1}\varphi_{(0)}(\xi)\mathcal{F}u\in M^p_q$ and $\mathcal{F}^{-1}\varphi_{j}(\xi)\mathcal{F}u\in M^p_q$ for every positive integer $j$, and that
\begin{equation*}
\|u|N^s_{p,q,r}\|:=\|\mathcal{F}^{-1}\varphi_{(0)}(\xi)\mathcal{F}u| M^p_q\|+ \|\{2^{sj}\|\mathcal{F}^{-1}\varphi_{j}(\xi)\mathcal{F}u|M^p_q\|\}^\infty_{j=1}|\ell^r\|<\infty, 
\end{equation*}
where $\mathcal{F}$ denotes the Fourier transform on ${\bf R}^N$. 
\end{Definition}
For every $t>0$ and every $u\in \mathcal{S} '({\bf R}^N)$, put $S(t)u:=\mathcal{F}^{-1}\exp(-t|\xi|^\theta)\mathcal{F}u$.
We formulate a solution to problem~\eqref{eq:1.1} and \eqref{eq:1.2} . 
\begin{Definition}
Let $T>0$ and $\varphi \in N^s_{p,q,r}$ for some $s\in{\bf R}$, $1\le q\le p<\infty$ and $1\le r\le \infty$. We say that $u$ is a solution to problem~\eqref{eq:1.1} in ${\bf R}^N\times [0,T)$ if 
$$u\in BC({\bf R}^N \times (\tau, T))$$
for $\tau\in (0,T)$, and $u$ satisfies 
$$u(x,t) = [S(t)\varphi](x) + \int_0^t [S(t-\tau)|u(\cdot,\tau)|^{\gamma-1}u(\cdot,\tau)](x)\,d\tau $$
for $(x,t)\in {\bf R}^N\times (0,T)$. 

\end{Definition}
\begin{Definition}
Let $T>0$ and $\varphi \in N^s_{p,q,r}$ for some $s\in{\bf R}$, $1\le q\le p<\infty$ and $1\le r\le \infty$. We say that $u$ is a solution to problem~\eqref{eq:1.2} in ${\bf R}^N\times [0,T)$ if 
$$u\, , \nabla u \in BC({\bf R}^N \times (\tau, T))$$
for $\tau\in (0,T)$, and $u$ satisfies 
$$u(x,t) = [S(t)\varphi](x) + \int_0^t [S(t-\tau)|\nabla u (\cdot,\tau)|^\gamma](x)\,d\tau $$
for $(x,t)\in {\bf R}^N\times (0,T)$. 

\end{Definition}

We are ready to state the main results of this paper. 
\begin{Theorem}
\label{Theorem:1.1}
Let $\gamma>1$, $\gamma \leq q\leq p<\infty$,  $-\theta/\gamma<s<0$ and $s\geq N/p-\theta/(\gamma-1)$. Then there exist  $\delta>0$  and $M>0$ such that for every $\varphi(x)\in N^s_{p,q,\infty}$ satisfying 
\begin{equation}
\label{eq:1.3}
\limsup_{j\to \infty}2^{sj}\|\mathcal{F}^{-1}\varphi_j\mathcal{F}\varphi|M^p_q\|<\delta,  
\end{equation}
problem~\eqref{eq:1.1} possesses the unique solution $u(x,t)$ on ${\bf R}^N \times [0,T)$ for some $T>0$ with a bound  $\sup_{0<t\leq T} t^{-s/\theta}\|u(\cdot,t)~|M^p_q\|\le M$.
\end{Theorem}
\begin{Remark}
\label{Remark:1.1}
To see the relation of these results with previous studies,  we remark here that inhomogeneous Besov-Morrey spaces under the assumption of Theorem~\ref{Theorem:1.1} includes the following functions and function spaces. Let $p_0= N(\gamma -1)/\theta$. 
\begin{itemize}
\item Let $\gamma>1+\theta/N$ and take $p$ as $\max \{\gamma, ~p_0\}<p<p_0\gamma $, then by Proposition~\ref{Proposition:2.1} and Proposition~\ref{Proposition:2.2}, we have 
\begin{equation*}
|x|^{-\frac{\theta}{\gamma-1}} \in M^{p_0}_{p_0 \gamma /p, \infty} \subset N^0_{p_0, p_0\gamma /p, \infty} \subset N^{N/p-\theta/(\gamma-1)}_{p,\gamma, \infty}.  
\end{equation*}
Since the assumption of Theorem~\ref{Theorem:1.1} is satisfied with $N^{N/p-\theta/(\gamma-1)}_{p,\gamma, \infty}$ for above $p$, we see by \eqref{eq:1.3} that there exists $c>0$ such that, if 
\begin{equation*}
|\varphi (x)|\le c |x|^{-\frac{\theta}{\gamma-1}}, \quad x\in{\bf R}^N, \\
\end{equation*}
then probolem~\eqref{eq:1.1} possesses a local-in-time solution. This result is consistent with that of \cite{IKO2} and thus the condition \eqref{eq:1.3} is necessary.
\item Let $\gamma=1+\theta/N$. Then by Proposition~\ref{Proposition:2.1} and Proposition~\ref{Proposition:2.2},  for any $p>1$ we have 
\begin{equation*}
L^p = M^p_p \subset N^0_{p,p,\infty}\subset N^{-N/p +N/\gamma}_{\gamma, \gamma, \infty}. 
\end{equation*}
Since the assumption of Theorem~\ref{Theorem:1.1} is satisfied with $N^{-N/p +N/\gamma}_{\gamma, \gamma, \infty}$, we see that if $\varphi \in L^p$ with $p>1$, then probolem~\eqref{eq:1.1} possesses a local-in-time solution. Note that in the case of $\theta=2$, problem~\eqref{eq:1.1} is not well-posed in $L^1$ (See for example, \cite{BC, CZ}).
\item
Let $1<\gamma<1+\theta/N$.  Then by Proposition~\ref{Proposition:2.1} and Proposition~\ref{Proposition:2.2}, we have 
\begin{equation*}
\delta (x) \in M^1\subset N^0_{1,1,\infty}\subset N^{-N+N/\gamma}_{\gamma, \gamma, \infty}.  
\end{equation*}
Since the assumption of Theorem~\ref{Theorem:1.1} is satisfied with $N^{-N+N/\gamma}_{\gamma, \gamma, \infty}$, we see that if $\varphi $ is a Radon measure, then probolem~\eqref{eq:1.1} possesses a local-in-time solution, which is consistent with the result of \cite{IKO2}. 
Furthermore, since  
\begin{equation*}
\partial^{|\alpha|}  \delta (x) \in N^{-N+N/\gamma-|\alpha|}_{\gamma, \gamma, \infty},  
\end{equation*}
we see that probolem~\eqref{eq:1.1} possesses a local-in-time solution for $\varphi = \partial ^{[\theta]} \delta$ and $\gamma<\frac{N+\theta}{N+[\theta]}$ if $\theta$ is not an integer, and for $\varphi = \partial ^{\theta-1} \delta$ and $\gamma<\frac{N+\theta}{N+\theta-1}$ if $\theta$ is  an integer. 
\end{itemize}

\end{Remark}

\begin{Theorem}
\label{Theorem:1.2}
Let $1<\gamma<\theta$, $\gamma \leq q\leq p<\infty$, $p>N(\gamma-1)/(\theta-1)$, $1-\theta/\gamma<s<0$ and $s\geq N/p+(\gamma-\theta)/(\gamma-1)$.  Then there exist  $\delta>0$  and $M>0$ such that for every $\varphi(x)\in N^s_{p,q,\infty}$ satisfying 
$$\lim \sup_{j\to \infty}2^{sj}\|\mathcal{F}^{-1}\varphi_j\mathcal{F}\varphi|M^p_q\|<\delta,$$  problem~\eqref{eq:1.2} possesses the unique solution $u(x,t)$ on ${\bf R}^N \times [0,T)$ for some $T>0$  with a bound $\sup_{0<t\leq T} t^{-s/\theta}\|u(\cdot,t)~|M^p_q\|\le M$ and $\sup_{0<t\leq T} t^{-s/\theta+1/\theta}\|\nabla u(\cdot,t)~|M^p_q\|\le M$.
\end{Theorem}
\begin{Remark}
\label{Remark:1.2}
To see the relation of these results with previous studies, we remark here that inhomogeneous  Besov-Morrey spaces under the assumption of Theorem~\ref{Theorem:1.2} includes the following functions and function spaces. Let $p_{1}= N(\gamma -1)/(\theta-\gamma)$. 
\begin{itemize}
\item Let $(N+\theta)/(N+1)<\gamma<\theta$ and   take $p$ as $\max \{\gamma, ~p_1\}<p<p_1\gamma $, then by Proposition~\ref{Proposition:2.1} and Proposition~\ref{Proposition:2.2}, we have 
\begin{equation*}
|x|^{-\frac{\theta-\gamma}{\gamma-1}} \in M^{p_1}_{p_1 \gamma /p, \infty} \subset N^0_{p_1, p_1\gamma /p, \infty} \subset N^{N/p-(\theta-\gamma)/(\gamma-1)}_{p,\gamma, \infty}. 
\end{equation*}
Since the assumption of Theorem~\ref{Theorem:1.2} is satisfied with $N^{N/p-(\theta-\gamma)/(\gamma-1)}_{p,\gamma, \infty}$ for above $p$, we see that there exists $c>0$ such that, if 
\begin{equation*}
|\varphi (x)|\le c |x|^{-\frac{\theta-\gamma}{\gamma-1}}, \quad x\in{\bf R}^N, \\
\end{equation*}
then probolem~\eqref{eq:1.2} possesses a local-in-time solution. This result is consistent with that of \cite{IKO2}.

\item Let $\gamma=(N+\theta)/(N+1)$. Then by Proposition~\ref{Proposition:2.1} and Proposition~\ref{Proposition:2.2},  for any $p>1$ we have 
\begin{equation*}
L^p = M^p_p \subset N^0_{p,p,\infty}\subset N^{-N/p +N/\gamma}_{\gamma, \gamma, \infty}. 
\end{equation*}
Since the assumption of Theorem~\ref{Theorem:1.2} is satisfied with $N^{-N/p +N/\gamma}_{\gamma, \gamma, \infty}$, we see that if $\varphi \in L^p$ with $p>1$, then probolem~\eqref{eq:1.2} possesses a local-in-time solution. 

\item
Let $1< \gamma < (N+\theta)/(N+1)$. Then by Proposition~\ref{Proposition:2.1} and Proposition~\ref{Proposition:2.2}, we have 
\begin{equation*}
\delta (x) \in M^1\subset N^0_{1,1,\infty}\subset N^{-N+N/\gamma}_{\gamma, \gamma, \infty}.  
\end{equation*}
Since the assumption of Theorem~\ref{Theorem:1.2} is satisfied with $N^{-N+N/\gamma}_{\gamma, \gamma, \infty}$, we see that if $\varphi $ is a Radon measure, then probolem~\eqref{eq:1.2} possesses a local-in-time solution. This result is consistent with that of \cite{IKO2}.
Furthermore, since 
\begin{equation*}
\partial^{|\alpha|}  \delta (x) \in N^{-N+N/\gamma-|\alpha|}_{\gamma, \gamma, \infty},  
\end{equation*}
we see that probolem~\eqref{eq:1.2} possesses a local-in-time solution for $\varphi = \partial ^{[\theta]-1} \delta$ and $\gamma<\frac{N+\theta}{N+[\theta]}$ if $\theta$ is not an integer, and for $\varphi = \partial ^{\theta-2} \delta$ and $\gamma<\frac{N+\theta}{N+\theta-1}$ if $\theta$ is  an integer. 
\end{itemize}

\end{Remark}

We explain the idea of the proof of Theorem~\ref{Theorem:1.1} and Theorem~\ref{Theorem:1.2}.
Let $S(t)u:=\mathcal{F}^{-1}\exp(-t|\xi|^\theta)\mathcal{F}u$.
 By modifying the arguments in \cite{KY}, we first prove the heat kernel estimates of the fractional Laplacian in inhomogeneous Besov-Morrey spaces and obtain
 the estimate
\begin{equation*}
\|S(t)u|N^\sigma_{p,q,1}\|\leq C(1+t^{(s-\sigma)/\theta})\|u|N^s_{p,q,\infty}\|,
\end{equation*}
\begin{equation*}
\|\nabla S(t)u|N^\sigma_{p,q,1}\|\leq C(1+t^{(s-\sigma-1)/\theta})\|u|N^s_{p,q,\infty}\|,
\end{equation*}
for $t>0$ and $\sigma>s$. Here, one of the main difficulties comes from the non-smoothness of the function $\exp(-t|\xi|^\theta)$, see Lemma~\ref{Lemma:2.2} and Remark~\ref{Remark:2.1}. 
 
 Then we show that the approximate solutions converge in some Banach space based on the local Morrey spaces with a bound near $t=0$. 

The rest of this paper is organized as follows. In Sections~2,  we obtain the heat kernel estimates of the fractional Laplacian in inhomogeneous Besov-Morrey spaces.  In section~3, we prove Theorem \ref{Theorem:1.1}.  In section~4, we prove Theorem \ref{Theorem:1.2}.

\section{Preliminaries}

    In this section,  we recall some preliminary facts about Besov-Morrey spaces and give estimates of heat kernel of fractional Laplacian in these function spaces. 
    
\noindent The following two propositions collect basic facts about Morrey spaces and Besov-Morrey spaces.

\begin{Proposition}[{\cite[Theorem~2.5]{KY}}]
\label{Proposition:2.1}
Let $1\leq q\leq p<\infty$, $r \in[1,\infty]$ and $s\in{\bf R}$. Then the following embeddings are continuous: 
\begin{align}
\label{eq:2.1}
&N^s_{p,q,r}\subset B^{s-N/p}_{\infty, r},\\
\label{eq:2.2}
&N^s_{p,q,r} \subset N^{s-N(1-l)/p}_{p/l, q/l, r} \quad {\text for \, \, any }\quad l \in (0,1). 
\end{align}
\end{Proposition}

\begin{Proposition}[{\cite[Proposition 2.11]{KY}}]
\label{Proposition:2.2}
Let $1\leq q\leq p<\infty$. Then the following embeddings are continuous:  
\begin{align}
\label{eq:2.3}
& N^0_{p,q,1}\subset M^p_q \subset N^0_{p,q,\infty}, \\
& M^p\subset N^0_{p,1,\infty}. \nonumber 
\end{align}
\end{Proposition}

We modify the arguments in \cite[Theorem~2.9\,(2)]{KY} and prepare the following two lemmas  for the estimates of heat kernel of fractional Laplacian in inhomogeneous Besov-Morrey spaces. Here, we denote by $\lfloor{x}\rfloor$ the greatest integer less than or equal to $x\in {\bf R}$. 
\begin{Lemma}
\label{Lemma:2.1}
Let $m\in {\bf R}$, $1\le q\le p<\infty$ and $P(\xi)\in   C^{\lfloor N/2 \rfloor+1}({\bf R}^N\setminus \{0\})$. Assume that there is $A>0$ such that 
\begin{equation*}
\left|\frac{\partial ^\alpha P}{\partial \xi^\alpha}(\xi)\right|\le A|\xi|^{m-|\alpha|}
\end{equation*}
for all $\alpha \in ({\bf N} \cup \{0\})^N$ with $|\alpha|\le \lfloor N/2 \rfloor+1$ and for all $\xi \neq 0$. Then the multiplier operator $P(D)u:= \mathcal{F}^{-1} P(\xi)\mathcal{F}u$  satisfies the estimate 
\begin{equation*}
\left\|\mathcal{F}^{-1}\varphi_j \mathcal{F}(P(D)u)|M^p_q\right\|\le CA2^{mj}\left\|\mathcal{F}^{-1}\varphi_j \mathcal{F}u|M^p_q\right\|
\end{equation*}
for every positive integer $j$ and $u\in \mathcal{S}'({\bf R}^N)$ such that $\mathcal{F}^{-1}\varphi_{j}(\xi)\mathcal{F}u\in M^p_q$, where $C>0$ is a constant independent of $j$, $A$, and $u$. 
\end{Lemma}
{\bf Proof.}
Put $\Phi_j:=\varphi_{j-1} + \varphi_j +\varphi_{j+1}$ and $K(x):=\mathcal{F}^{-1}\Phi_j(\xi)P(\xi)$ for $j\in {\bf Z}$. Note that ${\rm supp}~ \varphi_j(\xi)\subset \{\xi \in {\bf R}^N; 2^{j+1}/3\le |\xi| \le 2^{j+1}\}$ and $\Phi_j\equiv 1$ on ${\rm supp}~ \varphi_j(\xi)$. 

Putting also $N_0:=\lfloor N/2 \rfloor+1$, we have 
\begin{equation*}
\begin{split}
&\|K|L^1({\bf R}^N)\| 
= \int_{|x|\le2^{-j}}|K(x)|\,dx + \int_{|x|\ge2^{-j}}|K(x)| \,dx\\
& \le \left(\int_{|x|\le2^{-j}}\,dx\right)^{1/2}\left(\int_{|x|\le2^{-j}}|K(x)|^2\,dx\right)^{1/2}+\\
& \left(\int_{|x|\ge2^{-j}}|x|^{-2N_0}\, dx\right)^{1/2}\left(\int_{|x|\ge2^{-j}}|x|^{2N_0}|K(x)|^2\, dx\right)^{1/2}\\
& \le C\left(2^{-Nj/2}\|K(x)|L^2({\bf R}^N)\|+ 2^{(N_0-N/2)j}\sum_{|\alpha|=N_0}\|x^\alpha K(x)|L^2({\bf R}^N)\|\right)\\
&= C\left(2^{-Nj/2}\|\Phi_j(\xi)P(\xi)|L^2({\bf R}^N)\|+ 2^{(N_0-N/2)j}\sum_{|\alpha|=N_0}\left\|\frac{\partial^{|\alpha|}}{\partial \xi ^\alpha} (\Phi_j(\xi)P(\xi))|L^2({\bf R}^N)\right\|\right)\\
& \le C(2^{-Nj/2}2^{(m+N/2)j}A+2^{(N_0-N/2)j}2^{(m-N_0+N/2)j}A)=C2^{mj}A
\end{split}
\end{equation*}
for some constant $C>0$, depending on $N$, $m$, $\|\zeta|BC^{N_0}(R)\|$, but not on $j$ and $A$. 

Since $\mathcal{F}^{-1}\varphi_j \mathcal{F}(P(D)u)=K*(\mathcal{F}^{-1}\varphi_j \mathcal{F}u)$, 
 we see by \cite[Lemma~1.8]{KY} that 
$$\|\mathcal{F}^{-1}\varphi_j \mathcal{F}(P(D)u)|M^p_q\|\le CA2^{mj}\|\mathcal{F}^{-1}\varphi_j \mathcal{F}u|M^p_q\|$$
 for every positive integer $j$, and the proof is complete. 
$\Box$

\begin{Lemma}
\label{Lemma:2.2}
Let $m>0$, $1\le q\le p<\infty$ and $P(\xi)\in   C^{\lfloor N/2 \rfloor+1}({\bf R}^N\setminus \{0\})$. Assume that there is $A>0$ such that 
\begin{equation*}
\left|\frac{\partial ^\alpha P}{\partial \xi^\alpha}(\xi)\right|\le A|\xi|^{m-|\alpha|}
\end{equation*}
for all $\alpha \in ({\bf N} \cup \{0\})^N$ with $|\alpha|\le \lfloor N/2 \rfloor+1$ and for all $\xi \in B(0,4)\setminus \{0\}$. Then the multiplier operator $P(D)u:= \mathcal{F}^{-1} P(\xi)\mathcal{F}u$  satisfies the estimate 
\begin{equation*}
\|\mathcal{F}^{-1}\varphi_{(0)} \mathcal{F}(P(D)u)|M^p_q\|\le CA\|\mathcal{F}^{-1}\varphi_{(0)} \mathcal{F}u|M^p_q\|
\end{equation*}
for every $u\in \mathcal{S}'({\bf R}^N)$ such that $\mathcal{F}^{-1}\varphi_{(0)}(\xi)\mathcal{F}u\in M^p_q$, where $C>0$ is a constant independent of $A$ and $u$.

\noindent{\bf Proof.}
Put $K_j (x):=\mathcal{F}^{-1}\varphi_j(\xi)P(\xi)$ and $\Phi_{(0)}:=\varphi_{(0)}+ \varphi_1$. In the same way as in Lemma~\ref{Lemma:2.1}, we have 
\begin{equation*}
\begin{split}
\|\mathcal{F}^{-1}\Phi_{(0)}(\xi)P(\xi)|L^1({\bf R}^N)\| 
& \le \sum_{j=-\infty} ^1 \|K_j|L^1({\bf R}^N)\| \\
& \le \sum_{j=-\infty} ^1  C2^{mj}A \le CA
\end{split}
\end{equation*}
with some constant $C>0$ independent of $A$. This implies in the same way as in Lemma~\ref{Lemma:2.1}
\begin{equation*}
\|\mathcal{F}^{-1}\varphi_{(0)} \mathcal{F}(P(D)u)|M^p_q\|\le CA\|\mathcal{F}^{-1}\varphi_{(0)} \mathcal{F}u|M^p_q\|, 
\end{equation*}
and the proof is complete. 
$\Box$

\begin{Remark}
\label{Remark:2.1}
Note that we do not assume the smoothness of $P(\xi)$ at $\xi=0$, which is useful for the estimates of the derivative of heat kernel of fractional Laplacian since $P(\xi)=\exp (-t|\xi|^\theta)$ is not smooth at $\xi=0$ in general. In this respect, we improved \cite[Theorem~2.9\,(2)]{KY} where the smoothness at $\xi=0$ is needed. 
\end{Remark}
\end{Lemma}
In the following theorem, we obtain estimates of heat kernel of fractional Laplacian in inhomogeneous Besov-Morrey spaces. 
\begin{Theorem}
\label{Theorem:2.1}
Let $s\leq \sigma$, $1\leq q\leq p <\infty$ and  $r\in[1,\infty]$.  Then there exists $C>0$ such that the estimate
\begin{equation}
\label{eq:2.4}
\|S(t)u|N^\sigma_{p,q,r}\|\leq C(1+t^{(s-\sigma)/\theta})\|u|N^s_{p,q,r}\| \quad \text{for}\quad t>0
\end{equation}
holds. Furthermore, if $s<\sigma$, the estimate
\begin{equation}
\label{eq:2.5}
\|S(t)u|N^\sigma_{p,q,1}\|\leq C(1+t^{(s-\sigma)/\theta})\|u|N^s_{p,q,\infty}\| \quad \text{for}\quad t>0
\end{equation}
holds. 
\end{Theorem}

\noindent{\bf Proof.}
By induction we see that for every $\alpha \in {\bf N}^N$ there exist homogeneous polynomials $P_{\alpha,k} (\xi)$ of degree $|\alpha|$ for $k=1, 2, \ldots,|\alpha|$ such that for $\xi \neq 0$
\begin{equation}
\label{eq:2.6}
\frac{\partial ^{|\alpha|}\exp (-t|\xi|^\theta)}{\partial \xi_\alpha} = \exp (-t|\xi|^\theta)|\xi|^{-2|\alpha|}\sum_{k=1}^{|\alpha|}P_{\alpha, k}(\xi)t^k |\xi|^{k\theta}. 
\end{equation}
We have for  $m=s-\sigma$
\begin{equation*}
\begin{split}
|\xi|^{-m+|\alpha|}\frac{\partial ^{|\alpha|}\exp (-t|\xi|^\theta)}{\partial \xi_\alpha}
& \le C t^{\frac{m}{\theta}}\exp (-t|\xi|^\theta)\sum_{k=1}^{|\alpha|}(t^{\frac{1}{\theta}} |\xi|)^{k\theta-m}\\
& \le C_\alpha t^{\frac{m}{\theta}}. 
 \end{split} 
\end{equation*}
This together with Lemma~\ref{Lemma:2.1} implies 
\begin{equation}
\label{eq:2.7}
\|\mathcal{F}^{-1}\varphi_{j}(\xi)\mathcal{F}(S(t)u)|M^p_q\|\le Ct^{\frac{m}{\theta}}2^{mj}\|\mathcal{F}^{-1}\varphi_{j}(\xi)\mathcal{F}u|M^p_q\|
\end{equation}
for every positive integer and every $t>0$. 
On the other hand, since
\begin{equation*}
\begin{split}
& \|\mathcal{F}^{-1}\varphi_{(0)}(\xi)\mathcal{F}(S(t)u)|M^p_q\|\\
& \le\|\mathcal{F}^{-1}\Phi_{(0)}*\mathcal{F}^{-1}\exp(-t|\xi|^\theta)|L^1({\bf R}^N)|\|\|\mathcal{F}^{-1}\varphi_{(0)}(\xi)\mathcal{F}u|M^p_q\| \\
& \le C\|\mathcal{F}^{-1}\varphi_{(0)}(\xi)\mathcal{F}u|M^p_q\|, 
\end{split}
\end{equation*}
where $\Phi_{(0)}$ is as in Lemma~\ref{Lemma:2.2}. This together with \eqref{eq:2.7} implies the inequality $\eqref{eq:2.4}$.  
The inequality \eqref{eq:2.5} follows exactly in the same way as in \cite[Theorem~3.1]{KY} from the inequality \eqref{eq:2.4}, and the proof is complete. 
$\Box$

In the following lemma, we obtain another estimate of the heat kernel of fractional Laplacian by using the smallness condition on the initial data. 
\begin{Lemma}
\label{Lemma:2.3}
Let $1\leq q\leq p<\infty$ and $s<\sigma$. Then there exists $A>0$ such that, for every $u\in N^s_{p,q,\infty}$ and every $B>0$, satisfying
\begin{equation*}
A~ \limsup_{j\to \infty}2^{sj}\|\mathcal{F}^{-1}\varphi_j\mathcal{F}u|M^p_q\|<B,
\end{equation*}
there exists $T>0$ such that 
\begin{equation*}
\sup_{0<t\leq T}t^{(\sigma -s)/\theta} \|S(t)u|N^\sigma_{p,q,1}\|<B.
\end{equation*}
\end{Lemma}
{\bf Proof.}
Let $C_0$ be a positive constant satisfying the estimate
\begin{equation*}
\|S(t)u|N^\sigma_{p,q,1}\|\le C_0(1+t^{(s-\sigma)/\theta})\|u|N^s_{p,q,\infty}\|, 
\end{equation*}
and put $C_1=\max\{1,2\|\mathcal{F}^{-1}\varphi_0|L^1({\bf R}^N)\|\}$ and $A=C_0C_1$. 

\noindent Take $\delta>0$ such that
\begin{equation*}
\limsup_{j\to \infty}2^{sj}\|\mathcal{F}^{-1}\varphi_j\mathcal{F}u|M^p_q\|<\delta<B/A,
\end{equation*} 
then for some $m\in {\bf N}$, the estimate 
\begin{equation*}
2^{sj}\|\mathcal{F}^{-1}\varphi_j \mathcal{F}u|M^p_q\|\le \delta<B/A  
\end{equation*}
holds for every $j\ge m$.
Put $u_1=\mathcal{F}^{-1}\varphi_{(0)}(2^{-m}\cdot)\mathcal{F}u$ and $u_2=u-u_1$. Since
\begin{equation*}
\begin{split}
& {\rm supp}~ \varphi_{(0)}(2^{-m}\xi)\subset \left\{\xi \in {\bf R}^N; |\xi| \le \frac{5}{3}2^{m}\right\},\\
& {\rm supp}~ \varphi_j(\xi)\subset \left\{\xi \in {\bf R}^N; \frac{2^{j+1}}{3}\le |\xi| \le 2^{j+1}\right\},\\
& \varphi_{(0)}(2^{-m}\xi) \equiv 1 \quad {\rm on} \quad \{\xi \in {\bf R}^N; |\xi| \le 3\cdot2^{m-1}\},
\end{split}
\end{equation*}
we have 
\begin{equation*}
\mathcal{F}^{-1}\varphi_j \mathcal{F}u_1=\left\{
\begin{aligned}
 \mathcal{F}^{-1}\varphi_j \mathcal{F}u \quad & \text{for }\quad j\le m-1, \\
 \mathcal{F}^{-1}(\varphi_{m-1}+ \varphi_m)\varphi_j \mathcal{F}u\quad & \text{for }\quad j=m,m+1, \\
 0\quad  & \text{for }\quad j\ge m+2, 
\end{aligned}
\right.
\end{equation*}
and 
\begin{equation*}
\mathcal{F}^{-1}\varphi_j \mathcal{F}u_2=\left\{
\begin{aligned}
 0 \quad & \text{for }\quad j\le m-1, \\
 \mathcal{F}^{-1}(\varphi_{m+1}+ \varphi_{m+2})\varphi_j \mathcal{F}u\quad & \text{for }\quad j=m,m+1, \\
 \mathcal{F}^{-1}\varphi_j \mathcal{F}u \quad & \text{for }\quad j\ge m+2. 
\end{aligned}
\right.
\end{equation*}
It follows from the the definition of the constant $C_1$ and the fact $\|\mathcal{F}^{-1}\varphi_j|L^1\|=\|\mathcal{F}^{-1}\varphi_0|L^1\|$  that $\|u_2|N^s_{p,q,\infty}\|\le C_1\delta$. 
Therefore, we have 
\begin{equation}
\label{eq:2.8}
\begin{split}
& t^{(\sigma -s)/\theta}\|S(t)u_2|N^\sigma _{p,q,1}\|\le C_0(1+t^{(\sigma -s)/\theta})\|u_2|N^s_{p,q,\infty}\| \\
& \le C_0C_1 \delta(1+T^{(\sigma -s)/\theta})= A\delta (1+T^{(\sigma -s)/\theta})<\frac{A\delta +B}{2}
\end{split}
\end{equation}
for every $t\in(0,T]$, by taking $T>0$ sufficiently small. On the other hand, since $u_1\in N^{(\sigma + s)/2}_{p,q,\infty}$, we have the estimate 
\begin{equation}
\label{eq:2.9}
\begin{split}
& t^{(\sigma -s)/\theta}\|S(t)u_1|N^\sigma _{p,q,1}\|\le C_0(t^{(\sigma -s)/2\theta)}+ t^{(\sigma -s)/\theta)})\|u_1|N^{(s+\sigma)/2}_{p,q,\infty}\| \\
& \le C_0T^{(\sigma -s)/2\theta}(1+T^{(\sigma -s)/2\theta})\|u_1|N^{(s+\sigma)/2}_{p,q,\infty}\|<\frac{B-A\delta }{2}
\end{split}
\end{equation}
for every $t\in(0,T]$, by taking $T>0$ sufficiently small. We obtain the conclusion from \eqref{eq:2.8} and \eqref{eq:2.9},  and the proof is complete. 
$\Box$

\section{Proof of Theorem~\ref{Theorem:1.1}.}
In this section, we prove Theorem~\ref{Theorem:1.1} by using Theorem~\ref{Theorem:2.1}.
Let $X_T$ denote the set of Lebesgue measurable functions $u(x,t)$ on ${\bf R}^N \times (0,T)$ such that 
\begin{equation*}
 \|u|X_T\|:=\sup_{0<t< T} t^{-s/\theta}\|u(\cdot,t)~|M^p_q\|<\infty. 
\end{equation*}
Set $u_0(x,t)=[S(t)\varphi](x)$. Define $u_n(x,t)$ $(n=1,2,\ldots)$ inductively by 
\begin{equation}
\label{eq:3.1}
u_n(x,t) := u_0(x,t) + \int_0^t [S(t-\tau)|u_{n-1}(\cdot,\tau)|^{\gamma-1}u_{n-1}(\cdot,\tau)](x)\,d\tau.
\end{equation}
We prepare the following three lemmas for the proof of Theorem~\ref{Theorem:1.1}. 
\begin{Lemma}
Let $\gamma>1$, $T\le 1$, $\gamma \leq q\leq p<\infty$,  $-\theta/\gamma<s<0$ and $s\geq N/p-\theta/(\gamma-1)$. Then there exists  $C_2>0$ independent of $T$ such that 

\begin{equation*}
\|u_{n+1}|X_T\|\leq \|u_0|X_T\|+C_2\|u_n|X_T\|^\gamma
\end{equation*}
for $n=0,1, \ldots$. 
\end{Lemma}
{\bf Proof.}
By \eqref{eq:2.2}, \eqref{eq:2.3}, \eqref{eq:2.5} and \eqref{eq:3.1}, we see that 
\begin{equation*}
\begin{split}
\|u_{n+1}(\cdot,t)-u_0(\cdot,t)|M^p_q\| &\leq C\|u_{n+1}(\cdot,t)-u_0(\cdot,t)|N^0_{p,q,1}\| \\
& \leq C\int_0^t\|S(t-\tau)|u_n(\cdot,\tau)|^{\gamma-1}u_n(\cdot,\tau)|N^0_{p,q,1}\|\,d\tau\\
& \leq C\int_0^t\|S(t-\tau)|u_n(\cdot,\tau)|^{\gamma-1}u_n(\cdot,\tau)|N^{N(\gamma-1)/p}_{p/\gamma,q/\gamma,1}\|\,d\tau\\
& \leq C\int_0^t \{1+(t-\tau)^{-N(\gamma-1)/p\theta}\}\||u_n(\cdot,\tau)|^\gamma|N^{0}_{p/\gamma,q/\gamma,\infty}\|\,d\tau\\
&\leq C\int_0^t (t-\tau)^{-N(\gamma-1)/p\theta}\||u_n(\cdot,\tau)|^\gamma|M^{p/\gamma}_{q/\gamma}\|\,d\tau\\
&\leq C\int_0^t (t-\tau)^{-N(\gamma-1)/p\theta}\|u_n(\cdot,\tau)|M^{p}_{q}\|^\gamma \,d\tau\\
&\leq C\|u_n|X_T\|^\gamma \int_0^t (t-\tau)^{-N(\gamma-1)/p\theta}\tau^{s\gamma/\theta}\,d\tau\\
&\leq C t^{-N(\gamma-1)/p\theta +s\gamma/\theta+1}\|u_n|X_T\|^\gamma.
\end{split}
\end{equation*}
Therefore, we have 
\begin{equation*}
\begin{split}
\|u_{n+1}-u_0|X_T\|
&\leq Ct^{1+(\gamma-1)(s/\theta-N/p\theta)}\|u_n|X_T\|^\gamma\\
& \leq C\|u_n|X_T\|^\gamma
\end{split}
\end{equation*}
for $ T\leq 1$, and the proof is complete. 
$\Box$

\begin{Lemma}
\label{Lemma:3.2}
Let $\gamma>1$, $\gamma \leq q\leq p<\infty$, $-\theta/\gamma<s<0$ and $s\geq N/p-\theta/(\gamma-1)$. Then there exists  $C_3>0$  such that for every $\varphi(x)\in N^s_{p,q,\infty}$ satisfying $\lim \sup_{j\to \infty}2^{sj}\|\mathcal{F}^{-1}\varphi_j\mathcal{F}\varphi|M^p_q\|<\delta$ for some $\delta>0$, we can choose a positive number $T\le 1$ so small that the inequality $\|u_0|X_T\|<C_3\delta$ holds. Furthermore, we can choose $\delta$ so small that  $\sup_{n}\|u_n|X_T\|\le M$ for some $M>0$. 
\end{Lemma}
{\bf Proof.}
By Lemma~\ref{Lemma:2.3} with $B=A\delta$, we can take  $T\le 1$ such that the estimate 
\begin{equation*}
\sup_{0<t\le T} t^{-s/\theta}\|u_0|N^0_{p,q,1}\|<A\delta
\end{equation*}
holds. This together with \eqref{eq:2.3} implies $\|u_0|X_T\|<C_3\delta$ for some constant $C_3>0$. 
For $\delta>0$ satisfying
\begin{equation*}
2^\gamma C_2 C_3^\gamma \delta^{\gamma-1}<1,
\end{equation*}
we see by induction that 
\begin{equation*}
\sup_{n}\|u_n|X_T\|\le 2C_3\delta =:M,
\end{equation*}
 and the proof is complete. 
$\Box$

\begin{Lemma}
\label{Lemma:3.3}
Let $\gamma>1$, $\gamma \leq q\leq p<\infty$, $-\theta/\gamma<s<0$ and $s\geq N/p-\theta/(\gamma-1)$. Suppose that $\delta $ and $T\leq 1$ are small enough so that the assertion of Lemma~\ref{Lemma:3.2} holds. Then there exists a positive constant $C$ independent of $T$ such that 

\begin{equation*}
\|u_{n+2}-u_{n+1}|X_T\|\leq CM^{\gamma-1}\|u_{n+1}-u_{n}|X_T\|
\end{equation*}
for $n=0,1, \ldots$. 
\end{Lemma}
{\bf Proof.}
By \eqref{eq:2.2}, \eqref{eq:2.3}, \eqref{eq:2.5} and \eqref{eq:3.1}, we see that 
\begin{equation*}
\begin{split}
& \|u_{n+2}(\cdot,t)-u_{n+1}(\cdot,t)|M^p_q\| \leq C\|u_{n+2}(\cdot,t)-u_{n+1}(\cdot,t)|N^0_{p,q,1}\| \\
& \leq C\int_0^t\|S(t-\tau)(|u_{n+1}(\cdot,\tau)|^{\gamma-1}u_{n+1}(\cdot,\tau)-|u_{n}(\cdot,\tau)|^{\gamma-1}u_{n}(\cdot,\tau))|N^0_{p,q,1}\|\,d\tau\\
& \leq C\int_0^t\|S(t-\tau)(|u_{n+1}(\cdot,\tau)|^{\gamma-1}u_{n+1}(\cdot,\tau)-|u_{n}(\cdot,\tau)|^{\gamma-1}u_{n}(\cdot,\tau))|N^{N(\gamma-1)/p}_{p/\gamma,q/\gamma,1}\|\,d\tau\\
& \leq C\int_0^t (t-\tau)^{-\frac{N(\gamma-1)}{p\theta}}\||u_{n+1}(\cdot,\tau)|^{\gamma-1}u_{n+1}(\cdot,\tau)-|u_{n}(\cdot,\tau)|^{\gamma-1}u_{n}(\cdot,\tau)|N^{0}_{p/\gamma,q/\gamma,\infty}\|\,d\tau\\
&\leq C\int_0^t (t-\tau)^{-\frac{N(\gamma-1)}{p\theta}}\||u_{n+1}(\cdot,\tau)|^{\gamma-1}u_{n+1}(\cdot,\tau)-|u_{n}(\cdot,\tau)|^{\gamma-1}u_{n}(\cdot,\tau)|M^{p/\gamma}_{q/\gamma}\|\,d\tau\\
&\leq C\int_0^t (t-\tau)^{-\frac{N(\gamma-1)}{p\theta}}\||u_{n+1}(\cdot,\tau)-u_{n}(\cdot,\tau)|(|u_{n+1}(\cdot,\tau)|^{\gamma-1}+|u_n(\cdot,\tau)|^{\gamma-1})|M^{p/\gamma}_{q/\gamma}\|\,d\tau\\
&\leq CM^{\gamma-1}\int_0^t (t-\tau)^{-\frac{N(\gamma-1)}{p\theta}}\|u_{n+1}(\cdot,\tau)-u_{n}(\cdot,\tau)|M^p_q\|\,d\tau\\
&\leq CM^{\gamma-1}\|u_{n+1}-u_{n}|X_T\|\int_0^t (t-\tau)^{-N(\gamma-1)/p\theta}  \tau^{s\gamma/\theta}\,d\tau\\
&\leq C M^{\gamma-1}t^{-N(\gamma-1)/p\theta +s\gamma/\theta+1}\|u_{n+1}-u_{n}|X_T\|.
\end{split}
\end{equation*}
We used here \cite[Lemma~1.4]{KY}. 
Therefore, we have 
\begin{equation*}
\begin{split}
\|u_{n+2}-u_{n+1}|X_T\|
&\leq CM^{\gamma-1}t^{1+(\gamma-1)(s/\theta-N/p\theta)}\|u_{n+1}-u_{n}|X_T\|\\
& \leq CM^{\gamma-1}\|u_{n+1}-u_{n}|X_T\|, 
\end{split}
\end{equation*} 
and the proof is complete. 
$\Box$

\noindent{\bf Proof of Theorem~\ref{Theorem:1.1}.}

Take $\delta$ and $T$ so small that 
\begin{equation*}
\|u_{n+2}-u_{n+1}|X_T\|\leq \frac{1}{2}\|u_{n+1}-u_{n}|X_T\|
\end{equation*}
for $n=0,1, \ldots$, and we see that $u_n(x,t)$ converges in $X_T$. Set $u(x,t)$ as a limit of $u_n(x,t)$ in $X_T$ and we see that 
\begin{equation}
\label{eq:3.2}
u(x,t) := [S(t)\varphi](x) + \int_0^t [S(t-\tau)|u(\cdot,\tau)|^{\gamma-1}u(\cdot,\tau)](x)\,d\tau .
\end{equation}

We next prove that $u(x,t)\in L^\infty([\varepsilon, T]\times {\bf R}^N)$ for every $\varepsilon>0$. 
Let $n$ be the smallest integer greater than $N\gamma /\theta p$. Then we can take an increasing sequence of positive numbers $\{p_j\}^n_{j=1}$ such that $p_1=p$, $N/p_{j+1}>N/p_j - \theta/\gamma$ for every $j=1,2,\cdots, n-1$ and $N/p_n<\theta/\gamma$. We also define $\{q_j\}^n_{j=1}$ and $\{s_j\}^n_{j=1}$ as $q_1=q$, $q_{j+1}=p_{j+1}q_j/p_j$, $s_1=s$ and $s_{j+1}=N/p_{j+1}-N/p_j$.

By the obtained result,  we see that the solution $u(x,t)$ belongs to the spaces 
\begin{equation*}
L^\infty\left (\left[\frac{\varepsilon}{2n},T\right], M^p_q\right)\subset L^\infty\left(\left[\frac{\varepsilon}{2n},T \right], N^0_{p_1, q_1, \infty}\right)\subset L^\infty \left( \left[\frac{\varepsilon}{2n},T\right], N^{s_2}_{p_2, q_2, \infty} \right). 
\end{equation*}
Since $\gamma \le q_2 \le p_2$, $-\gamma/\theta < s_2<0$ and $s_2\ge N/p_2-\theta/(\gamma -1)$, we can apply the obtained result to see $u(x,t)\in L^\infty\left (\left[\frac{2\varepsilon}{2n},T\right], M^{p_2}_{q_2}\right)$. In the same way, since
\begin{equation*}
L^\infty\left (\left[\frac{j\varepsilon}{2n},T\right], M^{p_j}_{q_j}\right)\subset L^\infty\left(\left[\frac{j\varepsilon}{2n},T \right], N^0_{p_j, q_j, \infty}\right)\subset L^\infty \left( \left[\frac{j\varepsilon}{2n},T\right], N^{s_j}_{p_{j+1}, q_{j+1}, \infty} \right), 
\end{equation*}
where  $\gamma \le q_{j+1} \le p_{j+1}$, $-\gamma/\theta < s_{j+1}<0$ and $s_{j+1}\ge N/p_{j+1}-\theta/(\gamma -1)$, we have $u(x,t)\in L^\infty\left (\left[\frac{(j+1)\varepsilon}{2n},T\right], M^{p_{j+1}}_{q_{j+1}}\right)$ for $j=1,2,\cdots, n-1$. Therefore, we have $u(x,t)\in L^\infty\left (\left[\frac{\varepsilon}{2},T\right], M^{p_n}_{q_n}\right)$, where $p_n> N\gamma /\theta $. It follows from \eqref{eq:2.1} that 
\begin{equation}
\label{eq:3.3}
\begin{split}
& \left \|\int_{\varepsilon/2}^t S(t-\tau)|u(\cdot,\tau)|^{\gamma-1}u(\cdot,\tau)\,d\tau|L^\infty \right\| \\
& \le C \int_{\varepsilon/2}^t \left \| S(t-\tau)|u(\cdot,\tau)|^{\gamma-1}u(\cdot,\tau)|B^0_{\infty, 1} \right\|\,d\tau \\
& \le C \int_{\varepsilon/2}^t \left \| S(t-\tau)|u(\cdot,\tau)|^{\gamma-1}u(\cdot,\tau)|N^{N\gamma/p_{n}}_{p_n/\gamma, q_n/\gamma, 1} \right\|\,d\tau \\
& \le C \int_{\varepsilon/2}^t  \left(1+(t-\tau)^{-N\gamma/\theta p_n}\right)\left \||u(\cdot,\tau)|^{\gamma-1}u(\cdot,\tau)|N^{0}_{p_n/\gamma, q_n/\gamma, \infty} \right\|\,d\tau \\
& \le C \int_{\varepsilon/2}^t  (t-\tau)^{-N\gamma/\theta p_n}\left \||u(\cdot,\tau)|^{\gamma}|M^{p_n/\gamma}_{q_n/\gamma} \right\|\,d\tau \\
& \le C \int_{\varepsilon/2}^t  (t-\tau)^{-N\gamma/\theta p_n}\left \|u(\cdot,\tau)|M^{p_n}_{q_n} \right\|^{\gamma}\,d\tau \\
& \le C \left(t-\frac{\varepsilon}{2}\right)^{1-N\gamma/\theta p_n}\sup_{\varepsilon/2 \le \tau \le t}\left \|u(\cdot,\tau)|M^{p_n}_{q_n} \right\|^{\gamma}\,d\tau \\
& \le C T^{1-N\gamma/\theta p_n}\sup_{\varepsilon/2 \le \tau \le t}\left \|u(\cdot,\tau)|M^{p_n}_{q_n} \right\|^{\gamma}\,d\tau <\infty 
\end{split}
\end{equation}
for $\varepsilon/2 \le t\le T\le 1$. On the other hand, we have 
\begin{equation}
\label{eq:3.4}
\begin{split}
& \|S(t-\varepsilon/2)u(\cdot, \varepsilon/2)|L^\infty\|\le C\|S(t-\varepsilon/2)u(\cdot,\varepsilon/2)|B^0_{\infty,1}\|\\
& \le C\|S(t-\varepsilon/2)u(\cdot,\varepsilon/2)|N^{N/p}_{p, q, 1}\|\\
& \le C \left(1+(t-\varepsilon/2)^{-N/\theta p}\right)\left\||u(\cdot,\varepsilon/2)||M^{p}_{q} \right\|\\
& \le C (\varepsilon/2)^{-N/\theta p}\left\||u(\cdot,\varepsilon/2)||M^{p}_{q} \right\|<\infty 
\end{split}
\end{equation}
for $\varepsilon \le t\le T\le 1$. Since 
\begin{equation*}
u(x,t)=\left[S\left(t-\frac{\varepsilon}{2}\right)u\left(\cdot,\frac{\varepsilon}{2}\right)\right](x) + \int_{\varepsilon/2}^t \left[S(t-s)|u(\cdot,\tau)|^{\gamma-1}u(\cdot,\tau)\right](x)\,d\tau, 
\end{equation*}
this together with \eqref{eq:3.3} and \eqref{eq:3.4} implies that $u(x,t)\in L^\infty([\varepsilon, T]\times {\bf R}^N)$ for every $\varepsilon>0$. 

Finally, we prove the uniqueness of the solution. Assume that $u^{(1)}(x,t)$ and $u^{(2)}(x,t)$ are solutions to \eqref{eq:3.2} satisfying $\sup_{0\le t\le T}{t^{-s/\theta}\|u^{(j)}(\cdot,t)|M^p_q\|}<\infty$. Let $\overline{u}=u^{(1)}-u^{(2)}$ and $h(t)=\|\overline{u}(\cdot,t)|M^p_q\|$. Then exactly in the same way as in the proof of Lemma~\ref{Lemma:3.3}, we have
\begin{equation*}
\sup_{0< t \le T} t^{-s/\theta} h(t) \le CM^{\gamma-1}\sup_{0< t \le T} t^{-s/\theta} h(t)\le \frac{1}{2}\sup_{0< t \le T} t^{-s/\theta} h(t).
\end{equation*} 
Therefore, we see that $\overline{u}\equiv 0$, and the proof is complete. 
$\Box$

\section{Proof of Theorem~\ref{Theorem:1.2}.}
In this section, we prove Theorem~\ref{Theorem:1.2}.
Let $T>0$ be small and consider the Banach space 
\begin{equation*}
Y_T:=\{u(x,t) ~ \text{on}~ (0,T)\times {\bf R}^N: \|u~|Y_T\|<\infty \}, 
\end{equation*}
where
\begin{equation*}
 \|u~|Y_T\|:=\sup_{0<t< T} \{t^{-s/\theta}\|u(\cdot,t)~|M^p_q\|+t^{(-s+1)/\theta}\|\nabla u(\cdot,t)~|M^p_q\|\}.
\end{equation*}
Set $u_0(x,t)=[S(t)\varphi](x)$. Define $u_n (x,t)$ $(n=1,2,\ldots)$ inductively by 
\begin{equation}
\label{eq:4.1}
u_n(x,t) := u_0(x,t) + \int_0^t [S(t-\tau)|\nabla u_{n-1}(\cdot,\tau)|^\gamma](x)\,d\tau.
\end{equation}

For every $t>0$ and every $u\in \mathcal{S} '$, put $S_j(t)u:=\mathcal{F}^{-1}(i\xi_j)\exp(-t|\xi|^\theta)\mathcal{F}u$ for $1\le j \le N$. As in Section~2, we prove the derivative estimate for ${S}(t)$ in the following theorem. 
\begin{Theorem}
Let $s\leq \sigma$, $1\leq q\leq p <\infty$ and  $r\in[1,\infty]$.  Then there exists $C>0$ such that the estimate
\begin{equation}
\label{eq:4.2}
\|S_j(t)u|N^\sigma_{p,q,r}\|\leq C(1+t^{(s-\sigma-1)/\theta})\|u|N^s_{p,q,r}\| \quad \text{for}\quad t>0
\end{equation}
holds. Furthermore, if $s<\sigma$, the estimate
\begin{equation}
\label{eq:4.3}
\|S_j(t)u|N^\sigma_{p,q,1}\|\leq C(1+t^{(s-\sigma-1)/\theta})\|u|N^s_{p,q,\infty}\| \quad \text{for}\quad t>0
\end{equation}
holds. 
\end{Theorem}
{\bf Proof.}
By \eqref{eq:2.6} we see that for every $\alpha \in {\bf N}^N$ there exists a homogeneous polynomial $P_{\alpha,k} (\xi)$ of degree $|\alpha|$ for $k=1, 2, \ldots,|\alpha|$ and $P_{\alpha-e_j,k} (\xi)$ of degree $|\alpha|-1$ for $k=1, 2, \ldots,|\alpha|-1$ such that for $\xi \neq 0$
\begin{equation*}
\begin{split}
\frac{\partial ^{|\alpha|}(i\xi_j)\exp (-t|\xi|^\theta)}{\partial \xi_\alpha} 
& = i\xi_j \exp (-t|\xi|^\theta)|\xi|^{-2|\alpha|}\sum_{k=1}^{|\alpha|}P_{\alpha, k}(\xi)t^k |\xi|^{k\theta} \\
& + i\alpha_j \exp (-t|\xi|^\theta)|\xi|^{-2|\alpha|+2}\sum_{k=1}^{|\alpha|-1}P_{\alpha-e_j, k}(\xi)t^k |\xi|^{k\theta}. 
\end{split}
\end{equation*}
We have for  $m=s-\sigma$
\begin{equation*}
\begin{split}
|\xi|^{-m+|\alpha|}\frac{\partial ^{|\alpha|}(i\xi_j)\exp (-t|\xi|^\theta)}{\partial \xi_\alpha}
& \le C t^{\frac{m-1}{\theta}}\exp (-t|\xi|^\theta)\sum_{k=1}^{|\alpha|}(t^{\frac{1}{\theta}} |\xi|)^{k\theta-m+1}\\
& \le C_\alpha t^{\frac{m-1}{\theta}}. 
 \end{split} 
\end{equation*}
This together with Lemma~\ref{Lemma:2.1} implies 
\begin{equation}
\label{eq:4.4}
\|\mathcal{F}^{-1}\varphi_{j}(\xi)\mathcal{F}(S_j(t)u)|M^p_q\|\le Ct^{\frac{m}{\theta}}2^{mj}\|\mathcal{F}^{-1}\varphi_{j}(\xi)\mathcal{F}u|M^p_q\|
\end{equation}
for every positive integer and every $t>0$. 
On the other hand, by \eqref{eq:2.6} we have 
\begin{equation*}
\left|\frac{\partial ^{|\alpha|}(i\xi_j)\exp (-t|\xi|^\theta)}{\partial \xi_\alpha}\right| \le C_\alpha |\xi|^{1-|\alpha|} . 
\end{equation*}
for every $\xi \in B(0,4)\setminus \{0\}$. This together with Lemma~\ref{Lemma:2.2} implies 
\begin{equation}
\label{eq:4.5}
\|\mathcal{F}^{-1}\varphi_{(0)}(\xi)\mathcal{F}(S(t)u)|M^p_q\|\le C\|\mathcal{F}^{-1}\varphi_{(0)}(\xi)\mathcal{F}u|M^p_q\|. 
\end{equation}
The inequality \eqref{eq:4.2} follows from \eqref{eq:4.4} and \eqref{eq:4.5}. 
The inequality \eqref{eq:4.3} follows exactly in the same way as in \cite[Theorem~3.1]{KY} from the inequality \ref{eq:4.2}, and the proof is complete. 
$\Box$

\begin{Lemma}
\label{Lemma:4.1}
Let $1\leq q\leq p<\infty$ and $s<\sigma$. Then there exists $A>0$ such that, for every $u\in N^s_{p,q,\infty}$ and every $B>0$, satisfying 
\begin{equation*}
A\lim \sup_{j\to \infty}2^{sj}\|\mathcal{F}^{-1}\varphi_j\mathcal{F}u|M^p_q\|<B,
\end{equation*}
there exists $T>0$ such that 
$$\sup_{0<t\leq T}t^{(\sigma -s+1)/\theta} \|S_j(t)u|N^\sigma_{p,q,1}\|<B.$$
\end{Lemma}
{\bf Proof.}
Let $C_0$ be a positive constant satisfying the estimate
\begin{equation*}
\|S_j(t)u|N^\sigma_{p,q,1}\|\le C_0(1+t^{(s-\sigma-1)/\theta})\|u|N^s_{p,q,\infty}\|, 
\end{equation*}
and put $C_1=\max\{1,2\|\mathcal{F}^{-1}\varphi_0|L^1({\bf R}^N)\|\}$ and  $A=C_0C_1$.  Then there exists $m\in {\bf N}$ such that the estimate $2^{sj}\|\mathcal{F}^{-1}\varphi_j \mathcal{F}u\|\le \delta<B/A$ holds for every $j\ge m$. Put $u_1=\mathcal{F}^{-1}\varphi_{(0)}(2^{-m}\xi)\mathcal{F}u$ and $u_2=u-u_1$. Take $\delta>0$, $m\in{\bf N}$, $u_1$ and $u_2$ as in the proof of Lemma~\ref{Lemma:2.3}. Then we have 
\begin{equation}
\label{eq:4.6}
\begin{split}
& t^{(\sigma -s+1)/\theta}\|S_j(t)u_2|N^\sigma _{p,q,1}\|\le C_0(1+t^{(\sigma -s+1)/\theta)}\|u_2|N^s_{p,q,\infty}\| \\
& \le C_0C_1 \delta(1+T^{(\sigma -s+1)/\theta})= A\delta (1+T^{(\sigma -s+1)/\theta})<\frac{A\delta +B}{2}
\end{split}
\end{equation}
for every $t\in(0,T]$, by taking $T>0$ sufficiently small. On the other hand, since $u_1\in N^{(\sigma + s)/2}_{p,q,\infty}$, we have the estimate 
\begin{equation}
\label{eq:4.7}
\begin{split}
& t^{(\sigma -s+1)/\theta}\|S_j(t)u_1|N^\sigma _{p,q,1}\|\le C(t^{(\sigma -s)/2\theta)}+ t^{(\sigma -s)/\theta)})\|u_1|N^{(s+\sigma)/2}_{p,q,\infty}\| \\
& \le CT^{(\sigma -s)/2\theta}(1+T^{(\sigma -s)/2\theta})\|u_1|N^{(s+\sigma)/2}_{p,q,\infty}\|<\frac{B-A\delta }{2}
\end{split}
\end{equation}
for every $t\in(0,T]$, by taking $T>0$ sufficiently small. We obtain the conclusion from \eqref{eq:4.6} and \eqref{eq:4.7}. The proof is complete. 
$\Box$

We prepare the following three lemmas for the proof of Theorem~\ref{Theorem:1.2}. 

\begin{Lemma}
\label{Lemma:4.2}
Let $1<\gamma<\theta$, $T\le 1$, $\gamma \leq q\leq p<\infty$, $p>N(\gamma-1)/(\theta-1)$, $1-\theta/\gamma<s<0$ and $s\geq N/p+(\gamma-\theta)/(\gamma-1)$. Then there exists  $C_2>0$ independent of $T$ such that 

\begin{equation*}
\|u_{n+1}|Y_T\|\leq \|u_0|Y_T\|+C_2\|u_n|Y_T\|^\gamma
\end{equation*}
for $n=0,1, \ldots$. 
\end{Lemma}
{\bf Proof.}
By \eqref{eq:2.2}, \eqref{eq:2.3}, \eqref{eq:2.5} and \eqref{eq:4.1} we see that
\begin{equation*}
\begin{split}
\|u_{n+1}(\cdot,t)-u_0(\cdot,t)|M^p_q\| &\leq C\|u_{n+1}(\cdot,t)-u_0(\cdot,t)|N^0_{p,q,1}\| \\
& \leq C\int_0^t\|S(t-\tau)|\nabla u_n(\cdot,\tau)|^{\gamma}|N^0_{p,q,1}\|\,d\tau\\
& \leq C\int_0^t\|S(t-\tau)|\nabla u_n(\cdot,\tau)|^{\gamma}|N^{N(\gamma-1)/p}_{p/\gamma,q/\gamma,1}\|\,d\tau\\
& \leq C\int_0^t \{1+(t-\tau)^{-N(\gamma-1)/p\theta}\}\||\nabla u_n(\cdot,\tau)|^\gamma|N^{0}_{p/\gamma,q/\gamma,\infty}\|\,d\tau\\
&\leq C\int_0^t (t-\tau)^{-N(\gamma-1)/p\theta}\||\nabla u_n(\cdot,\tau)|^\gamma|M^{p/\gamma}_{q/\gamma}\|\,d\tau\\
&\leq C\int_0^t (t-\tau)^{-N(\gamma-1)/p\theta}\|\nabla u_n(\cdot,\tau)|M^{p}_{q}\|^\gamma \,d\tau\\
&\leq C\|u_n|Y_T\|^\gamma \int_0^t (t-\tau)^{-N(\gamma-1)/p\theta}\tau^{(s-1)\gamma/\theta }\,d\tau\\
&\leq C t^{-N(\gamma-1)/p\theta +(s-1)\gamma/\theta+1}\|u_n|Y_T\|^\gamma. 
\end{split}
\end{equation*}
In the same way, by \eqref{eq:2.2}, \eqref{eq:2.3},  \eqref{eq:4.1} and \eqref{eq:4.3} we see that
\begin{equation*}
\begin{split}
& \|\partial _j u_{n+1}(\cdot,t)-\partial _ju_0(\cdot,t)|M^p_q\| 
 \leq C\int_0^t\|S_j(t-\tau)|\nabla u_n(\cdot,\tau)|^{\gamma}|N^0_{p,q,1}\|\,d\tau\\
& \leq C\int_0^t\|S_j(t-\tau)|\nabla u_n(\cdot,\tau)|^{\gamma}|N^{N(\gamma-1)/p}_{p/\gamma,q/\gamma,1}\|\,d\tau\\
& \leq C\int_0^t \{1+(t-\tau)^{-N(\gamma-1)/p\theta- 1/\theta}\}\||\nabla u_n(\cdot,\tau)|^\gamma|N^{0}_{p/\gamma,q/\gamma,\infty}\|\,d\tau\\
&\leq C\int_0^t (t-\tau)^{-N(\gamma-1)/p\theta - 1/\theta}\||\nabla u_n(\cdot,\tau)|^\gamma|M^{p/\gamma}_{q/\gamma}\|\,d\tau\\
&\leq C\int_0^t (t-\tau)^{-N(\gamma-1)/p\theta- 1/\theta}\|\nabla u_n(\cdot,\tau)|M^{p}_{q}\|^\gamma \,d\tau\\
&\leq C\|u_n|Y_T\|^\gamma \int_0^t (t-\tau)^{-N(\gamma-1)/p\theta- 1/\theta}\tau^{(s-1)\gamma/\theta }\,d\tau\\
&\leq C t^{-N(\gamma-1)/p\theta -1/\theta +(s-1)\gamma/\theta+1}\|u_n|Y_T\|^\gamma.  
\end{split}
\end{equation*}
Therefore, we have 
\begin{equation*}
\begin{split}
\|u_{n+1}-u_0|Y_T\|
&\leq Ct^{1+(\gamma-1)(s/\theta-N/p\theta)-\gamma/\theta}\|u_n|Y_T\|^\gamma\\
& \leq C\|u_n|Y_T\|^\gamma
\end{split}
\end{equation*}
for $ T\leq 1$, and the proof is complete. 
$\Box$

\begin{Lemma}
\label{Lemma:4.3}
Let $1<\gamma<\theta$, $T\le 1$, $\gamma \leq q\leq p<\infty$, $p>N(\gamma-1)/(\theta-1)$, $1-\theta/\gamma<s<0$ and $s\geq N/p+(\gamma-\theta)/(\gamma-1)$. Then there exists  $C_3>0$  such that for every $\varphi(x)\in N^s_{p,q,\infty}$ satisfying $\lim \sup_{j\to \infty}2^{sj}\|\mathcal{F}^{-1}\varphi_j\mathcal{F}\varphi|M^p_q\|<\delta$ for some $\delta>0$ we can choose a positive number $T\le 1$ so small that the inequality $\|u_0|Y_T\|<C_0\delta$ holds. Furthermore, we can choose $\delta$ so small that  $\|u_n|Y_T\|\le M$ for some $M>0$. 
\end{Lemma}
{\bf Proof.}
By Lemma~\ref{Lemma:4.1} with $B=A\delta$, we can take a positive number $T\le 1$ such that the estimate 
\begin{equation*}
\sup_{0<t\le T} (t^{-s/\theta}\|u_0|N^0_{p,q,1}\|+ t^{(-s+1)/\theta}\|\nabla u_0|N^0_{p,q,1}\|)<A\delta
\end{equation*}
holds. This together with \eqref{eq:2.3} implies $\|u_0|Y_T\|<C_3\delta$ for some constant $C_3>0$. 

\noindent For $\delta>0$ satisfying
\begin{equation*}
2^\gamma C_2 C_3^\gamma \delta^{\gamma-1}<1,
\end{equation*}
we see by induction that 
\begin{equation*}
\sup_{n}\|u_n|X_T\|\le 2C_3\delta =:M,
\end{equation*}
 and the proof is complete. 
$\Box$

\begin{Lemma}
\label{Lemma:4.4}
Let $1<\gamma<\theta$, $T\le 1$, $\gamma \leq q\leq p<\infty$, $p>N(\gamma-1)/(\theta-1)$, $1-\theta/\gamma<s<0$ and $s\geq N/p+(\gamma-\theta)/(\gamma-1)$. Suppose that $\delta $ and $T\leq 1$ are small enough so that the assertion of Lemma~\ref{Lemma:4.3} holds. Then there exists a positive constant $C$ independent of $T$ such that 

\begin{equation*}
\|u_{n+2}-u_{n+1}|Y_T\|\leq CM^{\gamma-1}\|u_{n+1}-u_{n}|Y_T\|
\end{equation*}
for $n=0,1, \ldots$. 
\end{Lemma}
{\bf Proof.}
By \eqref{eq:2.2}, \eqref{eq:2.3}, \eqref{eq:2.5} and \eqref{eq:4.1} we see that
 
\begin{equation*}
\begin{split}
& \|u_{n+2}(\cdot,t)-u_{n+1}(\cdot,t)|M^p_q\| \leq C\|u_{n+2}(\cdot,t)-u_{n+1}(\cdot,t)|N^0_{p,q,1}\| \\
& \leq C\int_0^t\|S(t-\tau)(|\nabla u_{n+1}(\cdot,\tau)|^{\gamma}-|\nabla u_{n}(\cdot,\tau)|^{\gamma})|N^0_{p,q,1}\|\,d\tau\\
& \leq C\int_0^t\|S(t-\tau)(|\nabla u_{n+1}(\cdot,\tau)|^{\gamma}-|\nabla u_{n}(\cdot,\tau)|^{\gamma})|N^{N(\gamma-1)/p}_{p/\gamma,q/\gamma,1}\|\,d\tau\\
& \leq C\int_0^t \{1+(t-\tau)^{-N(\gamma-1)/p\theta}\}\||\nabla u_{n+1}(\cdot,\tau)|^{\gamma}-|\nabla u_{n}(\cdot,\tau)|^{\gamma}|N^{0}_{p/\gamma,q/\gamma,\infty}\|\,d\tau\\
&\leq C\int_0^t (t-\tau)^{-N(\gamma-1)/p\theta }\||\nabla u_{n+1}(\cdot,\tau)|^{\gamma}-|\nabla u_{n}(\cdot,\tau)|^{\gamma}|M^{p/\gamma}_{q/\gamma}\|\,d\tau\\
&\leq CM^{\gamma-1}\|u_{n+1}-u_{n}|Y_T\|\int_0^t (t-\tau)^{-N(\gamma-1)/p\theta}  \tau^{(s-1)\gamma/\theta}\,d\tau\\
&\leq C M^{\gamma-1}t^{-N(\gamma-1)/p\theta +(s-1)\gamma/\theta+1}\|u_{n+1}-u_{n}|Y_T\|.
\end{split}
\end{equation*}
In the same way, by \eqref{eq:2.2}, \eqref{eq:2.3}, \eqref{eq:4.1} and \eqref{eq:4.3} we see that
\begin{equation*}
\begin{split}
& \|\partial _j (u_{n+2}(\cdot,t)-u_{n+1}(\cdot,t))|M^p_q\| \leq C\|\partial _j (u_{n+2}(\cdot,t)-u_{n+1}(\cdot,t))|N^0_{p,q,1}\| \\
& \leq C\int_0^t\|S_j(t-\tau)(|\nabla u_{n+1}(\cdot,\tau)|^{\gamma}-|\nabla u_{n}(\cdot,\tau)|^{\gamma})|N^0_{p,q,1}\|\,d\tau\\
& \leq C\int_0^t\|S_j(t-\tau)(|\nabla u_{n+1}(\cdot,\tau)|^{\gamma}-|\nabla u_{n}(\cdot,\tau)|^{\gamma})|N^{N(\gamma-1)/p}_{p/\gamma,q/\gamma,1}\|\,d\tau\\
& \leq C\int_0^t \{1+(t-\tau)^{-N(\gamma-1)/p\theta-1/\theta}\}\||\nabla u_{n+1}(\cdot,\tau)|^{\gamma}-|\nabla u_{n}(\cdot,\tau)|^{\gamma}|N^{0}_{p/\gamma,q/\gamma,\infty}\|\,d\tau\\
&\leq C\int_0^t (t-\tau)^{-N(\gamma-1)/p\theta-1/\theta}\||\nabla u_{n+1}(\cdot,\tau)|^{\gamma}-|\nabla u_{n}(\cdot,\tau)|^{\gamma}|M^{p/\gamma}_{q/\gamma}\|\,d\tau\\
&\leq CM^{\gamma-1}\|u_{n+1}-u_{n}|Y_T\|\int_0^t (t-\tau)^{-N(\gamma-1)/p\theta-1/\theta}  \tau^{(s-1)\gamma/\theta}\,d\tau\\
&\leq C M^{\gamma-1}t^{-N(\gamma-1)/p\theta-1/\theta +(s-1)\gamma/\theta+1}\|u_{n+1}-u_{n}|Y_T\|.
\end{split}
\end{equation*}
We used here \cite[Lemma~1.4]{KY}. 
Therefore, we have 
\begin{equation*}
\begin{split}
\|u_{n+2}-u_{n+1}|Y_T\|
&\leq CM^{\gamma-1}t^{1+(\gamma-1)(s/\theta-N/p\theta)-\gamma/\theta}\|u_{n+1}-u_{n}|Y_T\|\\
& \leq CM^{\gamma-1}\|u_{n+1}-u_{n}|Y_T\|, 
\end{split}
\end{equation*} 
and the proof is complete. 
$\Box$

\noindent{\bf Proof of Theorem~\ref{Theorem:1.2}.}
Take $\delta$ and $T$ so small that 
\begin{equation*}
\|u_{n+2}-u_{n+1}|Y_T\|\leq \frac{1}{2}\|u_{n+1}-u_{n}|Y_T\|
\end{equation*}
for $n=0,1, \ldots$, and we see that $u_n(x,t)$ converges in $Y_T$. Set $u(x,t)$ as a limit of $u_n(x,t)$ in $Y_T$ and we see that 
\begin{equation}
\label{eq:4.8}
u(x,t) := u_0(x,t) + \int_0^t [S(t-\tau)|\nabla u(\cdot,\tau)|^{\gamma}](x)\,d\tau .
\end{equation}

We next prove that $u(x,t)\in L^\infty([\varepsilon, T]\times {\bf R}^N)$ and $\nabla u(x,t)\in L^\infty([\varepsilon, T]\times {\bf R}^N)$ for every $\varepsilon>0$. 
Let $n$ be the smallest integer greater than $N\gamma /(\theta-\gamma) p$. Then we can take an increasing sequence of positive numbers $\{p_j\}^n_{j=1}$ such that $p_1=p$, $N/p_{j+1}>N/p_j - (\theta-\gamma)/\gamma$ for every $j=1,2,\cdots, n-1$ and $N/p_n<(\theta -\gamma)/\gamma$. We also define $\{q_j\}^n_{j=1}$ and $\{s_j\}^n_{j=1}$ as $q_1=q$, $q_{j+1}=p_{j+1}q_j/p_j$, $s_1=s$ and $s_{j+1}=N/p_{j+1}-N/p_j$.

By the obtained result,  we see that the solution $u(x,t)$ and $\nabla u(x,t)$ belong to the spaces 
\begin{equation*}
L^\infty\left (\left[\frac{\varepsilon}{2n},T\right], M^p_q\right)\subset L^\infty\left(\left[\frac{\varepsilon}{2n},T \right], N^0_{p_1, q_1, \infty}\right)\subset L^\infty \left( \left[\frac{\varepsilon}{2n},T\right], N^{s_2}_{p_2, q_2, \infty} \right). 
\end{equation*}
Since $\gamma \le q_2 \le p_2$, $p_2>N(\gamma-1)/(\theta-1)$, $1-\gamma/\theta < s_2<0$ and $s_2\ge N/p_2-(\theta-\gamma)/(\gamma -1)$, we can apply the obtained result to see $u(x,t)\in L^\infty\left (\left[\frac{2\varepsilon}{2n},T\right], M^{p_2}_{q_2}\right)$ and $\nabla u(x,t)\in L^\infty\left (\left[\frac{2\varepsilon}{2n},T\right], M^{p_2}_{q_2}\right)$. In the same way, since
\begin{equation*}
L^\infty\left (\left[\frac{j\varepsilon}{2n},T\right], M^{p_j}_{q_j}\right)\subset L^\infty\left(\left[\frac{j\varepsilon}{2n},T \right], N^0_{p_j, q_j, \infty}\right)\subset L^\infty \left( \left[\frac{j\varepsilon}{2n},T\right], N^{s_j}_{p_{j+1}, q_{j+1}, \infty} \right), 
\end{equation*}
where  $\gamma \le q_{j+1} \le p_{j+1}$, $p_{j+1}>N(\gamma-1)/(\theta-1)$, $1-\gamma/\theta < s_{j+1}<0$ and $s_{j+1}\ge N/p_{j+1}-(\theta-\gamma)/(\gamma -1)$, we have $u(x,t)\in L^\infty\left (\left[\frac{(j+1)\varepsilon}{2n},T\right], M^{p_{j+1}}_{q_{j+1}}\right)$ for $j=1,2,\cdots, n-1$. Therefore, we have $u(x,t)\in L^\infty\left (\left[\frac{\varepsilon}{2},T\right], M^{p_n}_{q_n}\right)$ and $\nabla u(x,t)\in L^\infty\left (\left[\frac{\varepsilon}{2},T\right], M^{p_n}_{q_n}\right)$, where $p_n> N\gamma /(\theta-\gamma)$. It follows that 
\begin{equation}
\label{eq:4.9}
\begin{split}
& \left \|\int_{\varepsilon/2}^t S(t-\tau)|\nabla u(\cdot,\tau)|^{\gamma}\,d\tau|L^\infty \right\| \\
& \le C \int_{\varepsilon/2}^t \left \| S(t-\tau)|\nabla u(\cdot,\tau)|^{\gamma}|B^0_{\infty, 1} \right\|\,d\tau \\
& \le C \int_{\varepsilon/2}^t \left \| S(t-\tau)|\nabla u(\cdot,\tau)|^{\gamma}|N^{N\gamma/p_{n}}_{p_n/\gamma, q_n/\gamma, 1} \right\|\,d\tau \\
& \le C \int_{\varepsilon/2}^t  \left(1+(t-\tau)^{-N\gamma/\theta p_n}\right)\left \||\nabla u(\cdot,\tau)|^{\gamma}|N^{0}_{p_n/\gamma, q_n/\gamma, \infty} \right\|\,d\tau \\
& \le C \int_{\varepsilon/2}^t  (t-\tau)^{-N\gamma/\theta p_n}\left \||\nabla u(\cdot,\tau)|^{\gamma}|M^{p_n/\gamma}_{q_n/\gamma} \right\|\,d\tau \\
& \le C \int_{\varepsilon/2}^t  (t-\tau)^{-N\gamma/\theta p_n}\left \|\nabla u(\cdot,\tau)|M^{p_n}_{q_n} \right\|^{\gamma}\,d\tau \\
& \le C \left(t-\frac{\varepsilon}{2}\right)^{1-N\gamma/\theta p_n}\sup_{\varepsilon/2 \le \tau \le t}\left \|\nabla u(\cdot,\tau)|M^{p_n}_{q_n} \right\|^{\gamma}\,d\tau \\
& \le C T^{1-N\gamma/\theta p_n}\sup_{\varepsilon/2 \le \tau \le t}\left \|\nabla u(\cdot,\tau)|M^{p_n}_{q_n} \right\|^{\gamma}\,d\tau <\infty 
\end{split}
\end{equation}
for $\varepsilon/2 \le t\le T\le 1$. On the other hand, we have 
\begin{equation}
\label{eq:4.10}
\begin{split}
& \|S(t-\varepsilon/2)u(\cdot,\varepsilon/2)|L^\infty\|\le C\|S(t-\varepsilon/2)u(\cdot,\varepsilon/2)|B^0_{\infty,1}\|\\
& \le C\|S(t-\varepsilon/2)u(\cdot,\varepsilon/2)|N^{N/p}_{p, q, 1}\|\\
& \le C \left(1+(t-\varepsilon/2)^{-N/\theta p}\right)\left\||u(\cdot,\varepsilon/2)||M^{p}_{q} \right\|\\
& \le C (\varepsilon/2)^{-N/\theta p}\left\||u(\cdot,\varepsilon/2)||M^{p}_{q} \right\|<\infty 
\end{split}
\end{equation}
for $\varepsilon \le t\le T\le 1$. Since 
\begin{equation*}
u(x,t)=\left[S\left(t-\frac{\varepsilon}{2}\right)u\left(\cdot,\frac{\varepsilon}{2}\right)\right](x) + \int_{\varepsilon/2}^t \left[S(t-\tau)|\nabla u(\cdot,\tau)|^{\gamma}\right](x)\,d\tau, 
\end{equation*}
this together with \eqref{eq:4.9} and \eqref{eq:4.10} implies that $u(x,t)\in L^\infty([\varepsilon, T]\times {\bf R}^N)$ for every $\varepsilon>0$. 
We next prove that $\nabla u(x,t)\in L^\infty([\varepsilon, T]\times {\bf R}^N)$for every $\varepsilon>0$. 
It follows that 
\begin{equation}
\label{eq:4.11}
\begin{split}
& \left \|\int_{\varepsilon/2}^t S_j(t-\tau)|\nabla u(\cdot,\tau)|^{\gamma}\,d\tau|L^\infty \right\| \\
& \le C \int_{\varepsilon/2}^t \left \| S_j(t-\tau)|\nabla u(\cdot,\tau)|^{\gamma}|B^0_{\infty, 1} \right\|\,d\tau \\
& \le C \int_{\varepsilon/2}^t \left \| S_j(t-\tau)|\nabla u(\cdot,\tau)|^{\gamma}|N^{N\gamma/p_{n}}_{p_n/\gamma, q_n/\gamma, 1} \right\|\,d\tau \\
& \le C \int_{\varepsilon/2}^t  \left(1+(t-\tau)^{-N\gamma/\theta p_n-1/\theta}\right)\left \||\nabla u(\cdot,\tau)|^{\gamma}|N^{0}_{p_n/\gamma, q_n/\gamma, \infty} \right\|\,d\tau \\
& \le C \int_{\varepsilon/2}^t  (t-\tau)^{-N\gamma/\theta p_n-1/\theta}\left \||\nabla u(\cdot,\tau)|^{\gamma}|M^{p_n/\gamma}_{q_n/\gamma} \right\|\,d\tau \\
& \le C \int_{\varepsilon/2}^t  (t-\tau)^{-N\gamma/\theta p_n-1/\theta}\left \|\nabla u(\cdot,\tau)|M^{p_n}_{q_n} \right\|^{\gamma}\,d\tau \\
& \le C \left(t-\frac{\varepsilon}{2}\right)^{1-N\gamma/\theta p_n-1/\theta}\sup_{\varepsilon/2 \le \tau \le t}\left \|\nabla u(\cdot,\tau)|M^{p_n}_{q_n} \right\|^{\gamma}\,d\tau \\
& \le C T^{1-N\gamma/\theta p_n-1/\theta}\sup_{\varepsilon/2 \le \tau \le t}\left \|\nabla u(\cdot,\tau)|M^{p_n}_{q_n} \right\|^{\gamma}\,d\tau <\infty 
\end{split}
\end{equation}
for $\varepsilon/2 \le t\le T\le 1$. On the other hand, we have 
\begin{equation}
\label{eq:4.12}
\begin{split}
& \|S_j(t-\varepsilon/2)u(\cdot,\varepsilon/2)|L^\infty\|\le C\|S_j(t-\varepsilon/2)u(\cdot,\varepsilon/2)|B^0_{\infty,1}\|\\
& \le C\|S_j(t-\varepsilon/2)u(\cdot,\varepsilon/2)|N^{N/p}_{p, q, 1}\|\\
& \le C \left(1+(t-\varepsilon/2)^{-N/\theta p -1/\theta}\right)\left\||u(\cdot,\varepsilon/2)||M^{p}_{q} \right\|\\
& \le C (\varepsilon/2)^{-N/\theta p-1/\theta}\left\||u(\cdot,\varepsilon/2)||M^{p}_{q} \right\|<\infty 
\end{split}
\end{equation}
for $\varepsilon \le t\le T\le 1$. Since 
\begin{equation*}
\nabla u(x,t)=\left[S_j\left(t-\frac{\varepsilon}{2}\right)u\left(\cdot,\frac{\varepsilon}{2}\right)\right](x) + \int_{\varepsilon/2}^t \left[S_j(t-\tau)|\nabla u(\cdot,\tau)|^{\gamma}\right](x)\,d\tau, 
\end{equation*}
this together with \eqref{eq:4.11} and \eqref{eq:4.12} implies that $\nabla u(x,t)\in L^\infty([\varepsilon, T]\times {\bf R}^N)$ for every $\varepsilon>0$. 

Finally, we prove the uniqueness of the solution. Assume that $u^{(1)}(x,t)$ and $u^{(2)}(x,t)$ are solutions to \eqref{eq:4.8} satisfying 
$$\sup_{0\le t\le T}{t^{-s/\theta}\|u^{(j)}(\cdot,t)|M^p_q\|}+{t^{(-s+1)/\theta}\|\nabla u^{(j)}(\cdot,t)|M^p_q\|}<\infty.$$
 Let $\overline{u}=u^{(1)}-u^{(2)}$ and $h(t)=\|\overline{u}(\cdot,t)|M^p_q\|$. Then exactly in the same way as in the proof of Lemma~\ref{Lemma:4.4}, we have
\begin{equation*}
\begin{split}
\sup_{0< t \le T} \{t^{-s/\theta} h(t) +t^{(-s+1)/\theta} h(t)\} 
& \le CM^{\gamma-1}\sup_{0< t \le T} \{t^{-s/\theta} h(t) +t^{(-s+1)/\theta} h(t)\}\\
& \le \frac{1}{2}\sup_{0< t \le T} \{t^{-s/\theta} h(t) +t^{(-s+1)/\theta} h(t)\}.
\end{split}
\end{equation*} 
Therefore, we see that $\overline{u}\equiv 0$, and the proof is complete.  
$\Box$





\begin{thebibliography}{10}


\bibitem{AB}
L. Amour and M. Ben-Artzi, \textit{Global existence and decay for viscous Hamilton-Jacobi equations}, 
 Nonlinear Anal.,  31  (1998), 621--628.
 \bibitem{A}
 D. Andreucci, \textit{Degenerate parabolic equations with initial data measures}, 
 Trans. Amer. Math. Soc.,  349  (1997),  3911--3923.

 \bibitem{BP}
 P. Baras and M.  Pierre, \textit{Critre d'existence de solutions positives pour des \'{e}quations semi-lin\'{e}aires non monotones}, Ann. Inst. H. Poincar\'{e} Anal. Non Lin\'{e}aire,  2  (1985), 185--212.

 \bibitem{BSW}
 M. Ben-Artzi, P. Souplet, and F. B. Weissler, \textit{The local theory for viscous Hamilton-Jacobi equations in Lebesgue spaces}, 
 J. Math. Pures Appl.,   81  (2002),  343--378. 
 
 \bibitem{BC}
 H. Brezis and T. Cazenave, \textit{A nonlinear heat equation with singular initial data}, J. Anal. Math. 68 (1996), 277--304.
 
 \bibitem{CZ}
 C. Celik and Z. Zhou, \textit{No local $L^1$ solution for a nonlinear heat equation}, Comm. Partial Differential Equations 28 (2003), 1807--1831.
 
\bibitem{C2}
S. Cui, \textit{Local and global existence of solutions to semilinear parabolic initial value problems}, 
 Nonlinear Anal.,  43  (2001), 293--323.
\bibitem{DI}
J. Droniou and C. Imbert, \textit{Fractal first-order partial differential equations}, 
 Arch. Ration. Mech. Anal.,  182  (2006), 299--331.
\bibitem{FL} 
R. Filippucci and  S. Lombardi, \textit{Fujita type results for parabolic inequalities with gradient terms}, 
 J. Differential Equations,  268  (2020),  1873--1910.
 \bibitem{GP}
  V. A. Galaktionov and S. I. Pohozaev, \textit{Existence and blow-up for higher-order semilinear parabolic equations: majorizing order-preserving operators}, 
 Indiana Univ. Math. J.,  51  (2002), 1321--1338.
 \bibitem{GG}
 F. Gazzola and H.-C. Grunau, \textit{Global solutions for superlinear parabolic equations involving the biharmonic operator for initial data with optimal slow decay}, 
 Calc. Var. Partial Differential Equations,  30  (2007), 389--415.
 \bibitem{G}
 T. Ghoul, \textit{An extension of Dickstein's "small lambda'' theorem for finite time blowup}, 
 Nonlinear Anal.,  74  (2011), 6105--6115.
\bibitem{HI1}
K. Hisa and K. Ishige, \textit{Existence of solutions for a fractional semilinear parabolic equation with singular initial data}, 
 Nonlinear Anal.,  175  (2018), 108--132.
 \bibitem{HI2}
 K. Hisa and K. Ishige, \textit{Solvability of the heat equation with a nonlinear boundary condition}, 
 SIAM J. Math. Anal.,  51  (2019),  565--594.
 \bibitem{HIT}
 K. Hisa and K. Ishige,  J. Takahashi, \textit{Existence of solutions for an inhomogeneous fractional semilinear heat equation}, 
 Nonlinear Anal.,  199  (2020), 111920, 28.
 \bibitem{IKO1}
 K. Ishige and T.  Kawakami, and S.   Okabe, \textit{Existence of solutions for a higher-order semilinear parabolic equation with singular initial data}, 
 Ann. Inst. H. Poincar\'{e} Anal. Non Lin\'{e}aire,  37  (2020), 1185--1209.
\bibitem{IKO2}
K. Ishige, T. Kawakami, and S. Okabe, \textit{Existence of solutions to nonlinear parabolic equations via majorant integral kernel}, Nonlinear Anal, 223  (2022), 22 pp.
\bibitem{KW}
G. Karch and W.A. Woyczy\'{n}ski,  \textit{Fractal Hamilton-Jacobi-KPZ equations}, 
 Trans. Amer. Math. Soc.,  360  (2008), 2423--2442.
\bibitem{KY} 
H. Kozono and  M. Yamazaki, \textit{Semilinear heat equations and the Navier-Stokes equation with distributions in new function spaces as initial data}, 
 Comm. Partial Differential Equations,   19  (1994),  959--1014.
\bibitem{LN}
T.-Y. Lee and W.-M. Ni, \textit{Global existence, large time behavior and life span of solutions of a semilinear parabolic Cauchy problem}, 
 Trans. Amer. Math. Soc.,  333  (1992), 365--378.
 \bibitem{P}
 G. Ponce, \textit{Global existence of small solutions to a class of nonlinear evolution equations}, 
 Nonlinear Anal.,  9  (1985),  399--418.
\bibitem{QS}
P. Quittner and  P. Souplet, \textit{``Superlinear Parabolic Problems. Blow-up, Global Existence and Steady
States"}, Birkh\"auser Advanced Texts, Basel, 2007. 
\bibitem{R}
F. Ribaud, \textit{Semilinear parabolic equations with distributions as initial data}, 
 Discrete Contin. Dynam. Systems,  3  (1997), 305--316.
 \bibitem{Su}
 S. Sugitani, \textit{On nonexistence of global solutions for some nonlinear integral equations}, 
 Osaka Math. J.,  12  (1975), 45--51.
\bibitem{T}
 J. Takahashi, \textit{Solvability of a semilinear parabolic equation with measures as initial data}, Geometric properties for parabolic and elliptic PDE's, 
 257--276, Springer Proc. Math. Stat., 176, Springer,  2016. 
 \bibitem{W}
 F. B. Weissler, \textit{Existence and nonexistence of global solutions for a semilinear heat equation},
 Israel J. Math.,  38  (1981), 29--40.
 \bibitem{Wu}
 J. Wu,  \textit{Well-posedness of a semilinear heat equation with weak initial data}, 
 J. Fourier Anal. Appl.,  4  (1998), 629--642.
\end{thebibliography}
\end{document}